\documentclass[12pt]{article}
\usepackage[margin=1in,bmargin=1in,tmargin=1.2in]{geometry}
\usepackage{xcolor}
\usepackage{amsmath,amssymb,amsthm,nicematrix}
\usepackage{physics}
\usepackage{hyperref}
\usepackage{enumitem}
\usepackage{subcaption}
\usepackage{float}
\usepackage{parskip}
\usepackage{color,soul}
\usepackage{graphics}
\usepackage{pgf,tikz}
\usetikzlibrary{arrows,positioning} 
\def\flat{{\text{flat}}}
\newcommand{\ie}{\textit{i.e.}}
\newcommand{\BF}{\boldmath}
\newcommand{\mF}{F_{max}}

\newcommand{\R}{{\mathbb{R}}}

\newcommand{\NI}{{\noindent}}

\def\cF{{\cal F}}
\def\cL{{\cal L}}

\def\rank{\text{rank}}

\def\leb{\text{Leb}}
\def\cSx{{\text{SolSet}(x)}}
\def\cSp{{\text{SolSet}(p)}}
\tikzset{regnode/.style={draw,circle,scale=.8}}
\definecolor{ballblue}{rgb}{0.63, 0.79, 0.95}
\definecolor{bgreen}{rgb}{0.0, 0.7, 0.0}
\tikzset{rdnode/.style={draw,circle,color=red,fill=red!30,text=black}}
\tikzset{bdnode/.style={draw,circle,color=blue,fill=ballblue,text=black}}
\tikzset{bluenode/.style={draw,circle,color=blue,scale=.8,fill=ballblue,text=black}}
\tikzset{rednode/.style={draw,circle,scale=.8pt,color=red,fill=red!30,text=black}}
\tikzset{regnode/.style={draw,circle,scale=.8}}
\tikzset{uninode/.style={draw,circle,scale=.75pt}}
\tikzset{2uninode/.style={draw,circle,scale=.92pt}}
\tikzset{secnode/.style={draw,circle,scale=.75pt,color=red,fill=red!30,text=black}}
\tikzset{2secnode/.style={draw,circle,scale=.92pt,color=red,fill=red!30,text=black}}
\tikzset{tednode/.style={draw,circle,scale=.75pt,color=blue,fill=ballblue,text=black}}
\tikzset{2tednode/.style={draw,circle,scale=.92pt,color=blue,fill=ballblue,text=black}}
\tikzset{fnew/.style={draw,circle,scale=.75pt}}
\tikzset{fignew2/.style={draw,circle,scale=0.92pt}}
\theoremstyle{definition}
\newtheorem{theorem}{Theorem}
\newtheorem{definition}{Definition}
\newtheorem{lemma}{Lemma}
\newtheorem{proposition}[lemma]{Proposition}
\newtheorem{corollary}[lemma]{Corollary}

\newtheorem{example}{Example}
\newtheorem*{remark*}{Remark}
\title{ Robustness of solutions of almost every system of equations}
 \author{Sana Jahedi\footnote{Department of Mathematics and Statistics, University of New Brunswick. Corresponding author's email: s.jahedi@unb.ca }, Timothy Sauer\footnote{Department of Mathematical Sciences, George Mason University} and James A. Yorke\footnote{IPST, Mathematics, and Physics, University of Maryland College Park.}}
\date{}
\begin{document}
\maketitle

\begin{abstract} 
In mathematical modeling, it is common to have an equation $F(p)=c$ where the exact form of $F$ is not known. 
 This article shows that there are large classes of $F$
where almost all $F$ share the same properties. 
The classes we investigate are vector spaces $\mathcal{F}$ of $C^1$ functions $F:\mathbb{R}^N \to \mathbb{R}^M$ that satisfy the following condition:
$\mathcal{F}$ has ``almost constant rank'' (ACR) if 
there is a constant integer $\rho(\mathcal{F}) \geq 0$ such that rank$(DF(p))=\rho(\mathcal{F})$ for ``almost every'' $F\in \mathcal{F}$ and almost every $p\in\mathbb{R}^N$.
If the vector space $\mathcal{F}$ is finite-dimensional, then ``almost every'' is with respect to Lebesgue measure on $\mathcal{F}$, and otherwise, it means almost every in the sense of prevalence, as described herein.
Most function spaces commonly used for modeling purposes are ACR. 
In particular, we show that if all of the functions in $\mathcal{F}$ are linear or polynomial or real analytic,
or  if  $\mathcal{F}$ is the set of all functions in a ``structured system'', then $\mathcal{F}$ is ACR. For each $F$ and $p$, the solution set of $p\in\R^N$ is 
 $\cSp := \{x: F(x)=F(p)\}.$ A solution set of $F(p)=c$ is called robust if it persists despite small changes in $F$ and $c$. 
The following two global results are proved for almost every $F$ in an ACR vector space ${\cF}$: (1) Either the solution set $\cSp$ is robust for almost every $p\in\R^N$, or none of the solution sets are robust. (2) The solution set $\cSp$ is a $C^\infty$-manifold of dimension $d = N-\rho(\mathcal{F})$.
In particular, $d$ is the same for almost every $F\in\cF$.
\end{abstract}
\section{Introduction}
\begin{figure}[H]
    \centering
\begin{tikzpicture}[ultra thick]
  \coordinate (o) at (0,0);
     \node at (0,0) {$\begin{array}{r}
f_1(x_3)=c_1\\
f_2(x_3)=c_2\\
f_3(x_1,x_2,x_3)=c_3
\end{array}
  $};
\path ([yshift=-2ex]current bounding box.south) node[text width=4em] (sca){\subcaption{ }};
\node[regnode](s3)[xshift=13cm]{\bf 3};
\node[regnode](s1)[xshift=11cm]{\bf 1};
\node[regnode](s2)[xshift=15cm]{\bf 2};    
\draw[<->,>=stealth,black](s1)-- (s3);
\draw[<->,>=stealth,black](s2)-- (s3);
 \path node[text width=5em,right=7.3em of sca] (scb){\subcaption{ }};
\node[text width=5em,right=6.5em of scb] (scc){\subcaption{ }};
\draw[->,>=stealth,black](s3)edge[in=-20,out=60,loop above]node[below right]{}();
\node at (5,0){$DF = \left[
\begin{array}{cccc}
0& 0&f_{13}  \\
0 & 0&f_{23}\\
f_{31} & f_{32}& f_{33} 
\end{array}
\right]$};
    \end{tikzpicture}
    \caption{ \small
    {\bf 
    A fragile {structure} motivated by Competitive Exclusion Principle.} {\bf (a)} A structured system of equations describing the positive steady states of system of equations~\eqref{CE1}. For example, ${f_1}$ in the first equation is allowed to depend on ${x_3}$, but not ${x_1}$ or ${x_2}$. 
    This fact is represented in the two other parts of this figure. {\bf (b)} The structure matrix $DF(x) =\big[\frac{\partial f_i }{\partial x_j}(x)\big]$. {\bf (c)} The directed graph of the system. An edge from node $i$ to node $j$ in the graph means that variable $i$ is allowed to appear in equation $j$.
    Systems of this form cannot have a robust solution, so any solution that exists is fragile. 
    }
  \label{fg:Fragile}
\end{figure}
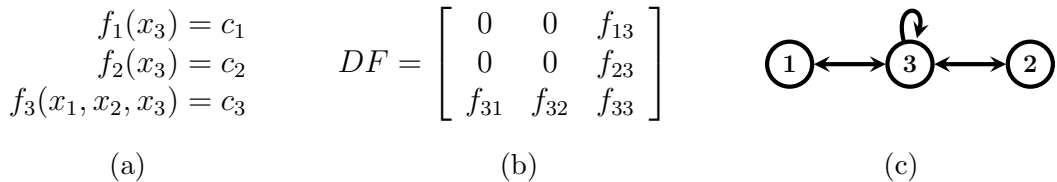

\NI The Competitive Exclusion Principle (CEP), long discussed in ecological literature,
holds that two predators that depend solely on the same prey species cannot coexist. More precisely, the CEP says the two predators whose population density depends purely on {the population density of a single prey species} cannot coexist unless they benefit precisely equally from the prey. {This scenario} is very unlikely in natural circumstances. These dynamics {can be described by the following system of equations.}
\begin{eqnarray} \label{CE1}
\frac{\dot{x}_1}{x_1}  &=&-c_1+f_1(x_3), \nonumber\\
\frac{\dot{x}_2}{x_2}  &=&-c_2+f_2(x_3),\\
\frac{\dot{x}_3}{x_3}  &=&-c_3+f_3(x_1,x_2,x_3),\nonumber
\end{eqnarray}
where $x_1$ and $x_2$ denote the population densities of the predators, and $x_3$ is the prey.
 To search for steady states of the system, set the left sides of the equations to zero, yielding the system of equations.
  \begin{eqnarray} \label{eqCEP}
 f_1(x_3) &=& c_1,\nonumber\\
 f_2(x_3) &=& c_2,\\
 f_3(x_1,x_2,x_3) &=& c_3,\nonumber
  \end{eqnarray}
illustrated in Fig.~\ref{fg:Fragile}(a). The first two equations share one unknown, $x_3$. 
 There may be a solution of these equations with positive $x_1, x_2, x_3$ for some exceptional $c=(c_1,c_2,c_3)$. However, it will fail to be robust in the sense that the set of solutions continues to exist under small perturbations of the equations without disappearing (or, more generally, changing the dimension of the set of solutions).
 In fact, for a dense subset of nearby choices of $c_1, c_2, f_1, f_2$, there will be no solutions.  
 
\NI The system of equations~\eqref{eqCEP} is a ``structured system'' of equations, in the sense that only certain variables are allowed to appear in certain equations. In this article, we prove some general facts about when robust solutions of structured systems can be expected and when they cannot, as in the above example. 

\NI {First note that if the system of equations~\eqref{eqCEP} has a solution such as $p=(p_1,p_2,p_3)$,} then $F(p) = (c_1,c_2,c_3)$. Thus we may write the system of equations~\eqref{eqCEP} as $F(x)=F(p)$ and consider the ``solution set'' for $p$,
\begin{align}\label{eq:solset}
    \cSp := \{x: F(x)=F(p)\}. 
\end{align}
For concreteness, we look at a particular example of system~\eqref{eqCEP}.
\begin{example} \rm
  \begin{eqnarray} \label{eqCEP1}
x_3^2  &=& c_1,  \nonumber\\
 x_3^4+1 &=& c_2, \\
 x_1^2-x_2+x_3^4 &=& c_3.\nonumber
  \end{eqnarray}
  Assume that $p$ is a solution, so $c_1 = p_3^2,  c_2 \ = \ p_3^4+1$, and $c_3 = p_1^2-p_2+p_3^4$.
Two facts are apparent: 

\NI (1) Although the system of equations~\eqref{eqCEP1} has a solution $p$, a dense subset of small $C^\infty$ perturbations of the equations $F$ will not have a solution. Thus this solution $p$ is not robust. We can see this from the first two equations alone, which would imply that $c_1^2+1 = c_2$, which only holds for special choices of $c$. 

\NI (2) If $p$ is a solution of the system of equations~\eqref{eqCEP1}, then SolSet$(p)$ is a one-dimensional $C^\infty$-manifold for every $p$. In fact, the first two equations will imply $x_3 = \pm p_3$. The last equation implies $x_1^2-x_2 = p_1^2-p_2$, Therefore 
\begin{align*}
 {\rm SolSet}(p) = \{(x_1,x_1^2-p_1^2+p_2,p_3)\} \cup \{(x_1,x_1^2-p_1^2+p_2,-p_3)\}.   
\end{align*}
\NI for all $x_1\in \mathbb{R}$,
which is a union of  one-dimensional curves (parabolas) in $R^3$.
\end{example}

\NI One of the main conclusions of this article is that the properties of Example 1 are quite general. Fix a structure, that pre-assigns particular variables to particular equations, as in (\ref{eqCEP}). Then for a dense subset of systems of $C^\infty$ functions  $F:U\subset R^N \to R^M$ with that structure, the following are true:

\begin{enumerate}
    \item Either 
\begin{itemize}
    \item[(a)] SolSet$(p)$ is robust for almost every $p$ in the domain $U$, or
    \item[(b)] SolSet$(p)$ is robust for no $p \in U$. See Corollary \ref{cor1}.
\end{itemize}
\item SolSet$(p)$ is a $C^\infty$-manifold for almost every $p\in U$.
See Theorem \ref{thm:flat}.
\end{enumerate}

\NI We will actually prove a stronger result, {that guarantees the above statements} for a ``prevalent'' subset of functions \cite{sauer_embedology_1991,hunt_prevalence_1992,ott_prevalence_2005}. Prevalent implies dense, and is a generalization of ``almost every'' to infinite dimensional function spaces. Furthermore, we show that distinguishing between cases 1(a) and 1(b) is a simple rank computation.
\NI A simple illustration of Case 1(b) is the example $F:\mathbb{R}^N\to\mathbb{R}$, $F(x)=x_1^2+x_2^2 + \ldots +  x_N^2=c$. Then $\cSp$ is an $(N-1)$-dimensional sphere for almost all $p\in\mathbb{R}^3.$ The only exception is $p=0$. Our generalization of this example will be the following:  For almost every real analytic function $F:\mathbb{R}^N\to\mathbb{R}$  and $p\in\mathbb{R}^N$,  rank$(DF(p)) = 1$, and  $\cSp$ is an $(N-1)$-manifold for almost every $F$ (in the sense of prevalence) and almost every $p$ (in the sense of Lebesgue measure).

The system of equations (\ref{eqCEP}) of CEP falls into category 1(b). The next example shows an instance from category 1(a).
\begin{example}
{\bf A robust example.} \rm The non-robust solution of  (\ref{eqCEP}) can be made robust.
Adding another prey species (node 4) to  Fig.~\ref{fg:Fragile}(c) yields the directed graph in Fig.~\ref{fg:rbstExpl}(c), and the equations
  \begin{eqnarray} \label{eqRob}
 f_1(x_3,x_4) &=& c_1,\nonumber\\
 f_2(x_3,x_4) &=& c_2,\\
 f_3(x_1,x_2,x_3) &=& c_3,\nonumber\\
 f_4(x_1,x_2) &=& c_4.\nonumber
  \end{eqnarray}
Let $\cF_2$ be the vector space of all $C^1$ functions $F=(f_1,f_2,f_3,f_4)$ where the $f_i$ are restricted to the form shown in Fig.~\ref{fg:rbstExpl}(a). 

We will show that for almost every $F\in\cF_2$ and for almost every $x=(x_1,\ldots,x_4)$ that is a solution of the system in Fig.~\ref{fg:rbstExpl}(a),  each sufficiently small perturbation of $c_1, c_2, c_3, c_4$  also has a solution. Thus, such solutions are allowed to exist in naturally-occurring circumstances.

\end{example}
\begin{figure}[H]
    \centering
     \begin{tikzpicture}[ultra thick]
  \coordinate (o) at (-1,0);
     \node at (-1,0) {$\begin{array}{r}
f_1(x_3, x_4)=c_1\\
f_2(x_3, x_4)=c_2 \\
f_3(x_1, x_2,x_3)=c_3\\
f_4(x_1, x_2)=c_4
\end{array}
  $};
   \path ([yshift=-2ex]current bounding box.south) node[text width=5em] (sca){\subcaption{ }};
 \node[regnode](s4)[xshift=12.5cm,yshift=-0.8cm]{\bf 4};
\node[regnode](s3)[xshift=12.4cm,yshift=0.7cm]{\bf 3};
\node[regnode](s1)[xshift=10.4cm,yshift=0.7cm]{\bf 1};
\node[regnode](s2)[xshift=14.5cm,yshift=0.7cm]{\bf 2};
\draw[->,>=stealth](s3)edge[in=-20,out=60,loop above]node[below right]{}();
\draw[<->,>=stealth](s2)-- (s4);
\draw[<->,>=stealth](s1)-- (s3);
\draw[<->,>=stealth](s2)-- (s3);
\draw[<->,>=stealth](s1)-- (s4);
  \path node[text width=5em,right=8em of sca] 
  (scb){\subcaption{ }};
\node[text width=5em,right=7.1em of scb] (scc){\subcaption{ }};
   \node at (4,0){$DF = \left[
\begin{array}{cccc}
0& 0 &f_{13} & f_{14}\\
0 & 0 & f_{23} & f_{24}\\
 f_{31}& f_{32} & f_{33} & 0\\
f_{41} &  f_{42} & 0 & 0
\end{array}
\right]$};
    \end{tikzpicture} 
    \caption{\small {\bf A robust family of systems}.
   For almost every $F$ of this form, $DF(x)$ is nonsingular, so the system will have robust solutions.  (a) The structured system of equations. (b) The Jacobian is generically of rank 4. {\bf (c)} The directed graph of the system.
}
  \label{fg:rbstExpl}
\end{figure}

\NI An even simpler case of Case 1(a)  is the real analytic equation $\cos x = c$, where $F:\mathbb{R}\to \mathbb{R}$ is given by $F(x) = \cos x$. For almost every $p$ in the domain, $F'(p)$ is nonzero, and the Implicit Function Theorem implies that the solution set is nonempty for small changes in $c$.

\NI In applications, there is  often uncertainty in the exact details of the equations.  
That motivates our focus on vector spaces of systems for each fixed structure, such as the structures in Figs.~\ref{fg:Fragile}  and \ref{fg:rbstExpl}. More generally, we will {propose a condition on} vector spaces $\mathcal{F}$
that implies that, for almost all $F\in\mathcal{F}$ and $p\in U$, the solution sets~\eqref{eq:solset} have the same properties,
such as almost all being robust, fragile, or having the same dimension.
When the vector space is infinite-dimensional, we will rely on the concept of prevalence to define ``almost every'', which is stronger than (implies) the notion of dense subset of the vector space of $C^\infty$ functions. A review of facts about prevalence is given in Appendix \ref{prevalence}.

\NI In the next section we define our terminology. In Section~\ref{Sec:applications} we offer some applications, in Section~\ref{s:proofs} we provide proofs of the theorems, and discuss relations with prior work in Section~\ref{s:discussion}.

\section{Function spaces of almost constant rank}

{Any} vector space $\cL$ of $M\times N$ matrices has the property that almost every matrix in the vector space has the same rank, equal to $\max_{A\in\cL} \rank(A)$, see Prop.~\ref{prop:matrices}.  We will argue that this key fact extends to many examples of vector spaces of  nonlinear functions.

{\begin{definition}
\normalfont
Let $F:U\subset R^N\to R^M$ be a $C^1$ function where $U$ is an open subset of $\mathbb{R}^N$.
We say  $F$ has {\bf almost constant rank (ACR)} if there is an integer {$\rho \geq 0$} such that rank$(DF(x)) = \rho$ for almost every $x\in U$. Define:
\begin{align}
 {\rm maxrank}(F) = \max_{x\in U}\rank\qty(DF(x)).   
\end{align}
If $F$ is an ACR function such that rank$\qty(DF(x))=\rho$ for almost every $x \in U$, then $\rho = $maxrank$(F)$, because $\rank(DF(x))$ takes its maximum value on an open set of $x$.
\end{definition}}

\NI 
We note below that real analytic functions have almost constant rank when $U$ is connected. However, there are $C^\infty$ functions that do not have almost constant rank, such as a monotonic $C^\infty$ function
 $F:\mathbb{R} \to\mathbb{R}$ for which $\frac{dF}{dx}(x)=0$ if and only if $x\le 0$. 
Then $\rank(DF(x)) =1$ if $x>0$ and 0 otherwise.

\begin{definition}\normalfont
For a vector space $\mathcal{F}$ of functions $F$, define: 
\begin{align}
{\rm maxrank}({\mathcal{F}}) =\max_{F\in\mathcal{F},  x\in U}{ \rm rank}(DF(x)).   
\end{align}
We say the vector space $\mathcal{F}$ has {\bf almost constant rank (ACR)} if for almost every $F\in\mathcal{F}$ and almost every $x$, ${ \rm rank}(DF(x)) =$ maxrank$({\mathcal{F}})$.
 \end{definition}
 
{\begin{definition}
\rm  We call an ACR function $F\in\mathcal{F}$  a 
 {\bf rank-maximizer} for the vector space $\mathcal{F}$ if maxrank$(F)=$ maxrank$(\mathcal{F})$.
\end{definition}}

 Throughout this paper, we deal with vector spaces $\mathcal{F}$ of functions that can be finite or infinite dimensional. By the phrase {\bf \BF ``almost every $ F\in\mathcal{F}$''},
when $\mathcal{F}$ is finite dimensional, we mean almost every with respect to the Lebesgue measure on the vector space $\mathcal{F}$. If $\mathcal{F}$ is infinite dimensional, we mean almost every in the sense of prevalence. See Appendix~\ref{prevalence} for a short primer on prevalence. 

 \NI It is clear that every ACR vector space contains a rank-maximizer. The next theorem states the converse. See Section \ref{s:proofs} for the proof.
\begin{theorem}
[ACR Vector Space]\label{thm:key}
 Let ${\mathcal{F}}$ be a vector space of $C^1$ functions
that has a rank-maximizer.
Then  ${\mathcal{F}}$ is ACR.
\end{theorem} 
\begin{definition}\normalfont
 Let $F:U\subseteq \mathbb{R}^N\to \mathbb{R}^M$.
We say a point $x$ is {\bf robust} for $F$
if $DF(x)$ has rank $M$, and $x$ is
{\bf fragile} if it is not robust.  
\end{definition}
This definition is suggested by the Implicit Function Theorem, which says that
if rank $DF(x) = M$, {each}
solution of $F(x) = c$ persists for small changes in $F$ and $c$.
The following result is an immediate consequence of the definition.
\begin{corollary}\label{cor1}  Let $F:U \subseteq \mathbb{R}^N \to \mathbb{R}^ M$ be an ACR function. Then
either 
\begin{itemize}
    \item[(a)] SolSet$(p)$ is robust for almost every $p$ in the domain $U$, or
    \item[(b)] SolSet$(p)$ is robust for no $p \in U$.
\end{itemize}
\end{corollary}

Many typical function spaces are ACR, including any vector space containing only real analytic functions.  
This fact follows from Theorem \ref{thm:key} and Proposition \ref{prop3}, which in turn depends on a basic fact about real analytic functions for which Mityagin published an elementary accessible proof in 2015 \& 2020. 
 \begin{proposition}\label{prop:polynomial}[See \cite{mityagin2015zero} ]. Let $U$ be an open connected subset of $\mathbb{R}^N$. Let $F:U\to\mathbb{R}^M$ be a real analytic function that is not identically zero. Then the set of solutions of $F(x_1,\ldots,x_N)=0$ has Lebesgue measure zero.
 \end{proposition}
\begin{proposition}\label{prop3}
Let $U$ be an open connected subset of $\mathbb{R}^N$. Let $F:U \to \mathbb{R}^M$ be a real analytic function. Then $\rank\qty(DF(x))$ is constant for almost every $x\in U$.  Hence, $F$ is ACR.
\end{proposition}
\begin{proof}
Let $\rho$ be the maximal rank of $DF(x)$ for $x\in \mathbb{R}^N$.
Hence, for some $x$,  $DF(x)$ has a $\rho\times\rho$ minor with determinant given by a real analytic function $P(x_1,\ldots,x_N)$ that is not identically zero. By Prop.\ref{prop:polynomial}, {the function $P(x_1,\ldots,x_N)$ is zero} only on a set of measure zero, and elsewhere it is nonzero. 
\end{proof} 
Since Prop.\ref{prop3} guarantees that
every vector space of real analytic functions on a connected open set has a rank-maximizer, it follows that any  vector space ${\mathcal{F}}$ consisting entirely of either linear functions or polynomial functions always has a rank-maximizer. In each of these cases, Theorem \ref{thm:key} implies the vector space is ACR.
In this sense such vector spaces are like spaces of linear maps.

{
\begin{definition} \rm

A {\bf structure matrix} S is a matrix where certain entries are allowed to be nonzero and the rest are zero.
A function F respects a structure matrix $S$ where $\frac{\partial F_i}{\partial x_j}(x) = 0$ for all $x$ when $S_{ij}=0$. In particular, let $\mathcal{L}(S)$ be the set of all functions $Ax$
where $A$ is a matrix that respects $S$.
{We say a vector space $\mathcal{F}$ of $C^1$ functions that respect a structure is a {\bf structured function space}, provided $\mathcal{F}$ includes $\mathcal{L}(S)$.
}
\end{definition}

Examples of structure matrices are found in Figures~\ref{fg:Fragile}(b)  and \ref{fg:rbstExpl}(b).}

\begin{proposition} \label{prop:str}
Every structured function space $\mathcal{F}$ is ACR.
In particular, the vector space of all $C^\infty$ functions that respect a structure is ACR.
\end{proposition}

\begin{proof}
{Choose $F$ in a structured function space $\mathcal{F}$ and $x\in U$
so that rank$(DF(x)) =$ maxrank$(\mathcal{F})$. Then
from the definition of 
structured function space, 
the matrix $A:= DF(x)$ also respects the structure and is in $\mathcal{F}$
 and its Jacobian $DA(x) = A$ is independent of $x$.
 Hence $A$ is a constant-rank function $A\in\mathcal{F}$ whose rank is maxrank$(\mathcal{F})$, so it is a rank maximizer.
Since $A$ has constant rank,} 
by Theorem \ref{thm:key}, $\mathcal{F}$ is ACR.
\end{proof}


\begin{example}\normalfont
 {\bf A vector space of functions that is not a structured system.} 
\normalfont
For fixed non-zero constants $a$ and $b$, let $\mathcal{F}$ be a vector space of functions $F(x):=(f_1,\ldots,f_4)$ of the form:\\
\begin{align}
f_1(x_1,x_2,x_3,x_4) &=c_1,\nonumber\\
f_2(x_1,x_2,x_3,x_4) &=c_2,\\
f_3(ax_1+bx_2) &=c_3,\nonumber\\
f_4(ax_1+bx_2) &=c_4.\nonumber
\end{align}\label{eq DF3}
 Imagine that this system represents two predators, species 3 and 4, and two prey species 1 and 2, with the assumption that the two prey species both provide the predators with the same nutrition which is proportional to $ax_1+bx_2$. So the predators are competing what becomes a single resource $ax_1+bx_2$ even if the prey look very different. One of the predators will almost certainly die, according to the Competitive Exclusion Principle.

For $F\in\mathcal{F}$, there are no robust solutions because rank$(DF(x))<4$ for all $x\in U$. 
{Since the functions $f_1$ and $f_2$ are functions of a single variable,  we write $f'_1$ and $f'_2$ for their derivatives.}
In fact, it is clear from the Jacobian

 \begin{equation}
 DF = \left[
\begin{array}{cccc}
f_{11}& f_{12}&f_{13}&f_{14} \\
f_{21}& f_{22}&f_{23}&f_{24} \\
a f'_{3} & bf'_{3}& 0&0\\
a f'_{4} & b f'_{4}& 0&0
\end{array}
\right]
 \end{equation}
 that maxrank$(F)$ is three for $F\in \mathcal{F}$. 
 {In particular $\left[
\begin{array}{cc}
a f'_{3} & bf'_{3}\\
a f'_{4} & b f'_{4}
\end{array}\right]$ has determinant zero.}
\NI For this example, we aim at simplicity. For nonzero $a$ and $b$, there is a change of variables that changes $\mathcal{F}$ into a structured system, but that is not true of more complicated systems.
\end{example}

{Theorem \ref{thm:flat} below and its corollary characterize the regularity of general solution sets, whether they are robust or fragile.} 
The following definition says a function $F$ essentially partitions the domain $U$ into manifolds when the sets $\cSp$ are manifolds for almost every $p$.

\begin{definition}\normalfont
{Let $U$ be an open subset of $\mathbb{R}^N$.} We say $F:U \to \mathbb{R}^M$ is {\bf\BF ``$d$-\flat''} if  $\cSp$ is a  $d$-dimensional manifold for  almost every $p\in U$. We say $F$ is ``\flat'' if it is
``d-\flat'' for some $d$.\end{definition}

This definition is also motivated by the Implicit Function Theorem. If rank $DF(x) = \rho$, then SolSet$(x)$ is an $(N- \rho)$-dimensional $C^\infty$ manifold in a neighborhood of $x$.
In this article, we use the term {\bf manifold} to mean $C^\infty$ manifold without boundary. 
Hence if $U=\mathbb{R}^N$, a sphere $\{x:\|x\|= 1\}$ is a manifold, but a closed disk $\{x:\|x\|\le 1\}$ is not.
A manifold can be either compact or unbounded. It can have a countable number of components, as in the $M=N=1$ case of $F(x)=\cos x$,  where ${\cSx}$ is a 0-manifold for each $x\in \mathbb{R}$  having countably many components. 

The archetypal example of a \flat\ function $F:\mathbb{R}^N\to\mathbb{R}^M$ is a linear function. 
 Let $\rho$ be the dimension of $F(\mathbb{R}^N)$.
For any $p\in \mathbb{R}^N$, ${\cSp}$ is a $(N-\rho)$-dimensional hyperplane. 
In the general definition of \flat, we allow a measure zero set of exceptional $p$.
Another example of a nonlinear \flat\ function is the equation $F(x_1,x_2) = x_1^2 + 2x_2^2$.
Here SolSet$(x_1,x_2)$ is an ellipse for all 
$(x_1,x_2)\ne (0,0)$.

As an example of a non-\flat\ $C^\infty$ function where $M=N=1$, consider 
$F(x)$ which is zero for $x \leq 0$, and 
 is strictly monotonically increasing for $x>0$, such as 
$F(x) =\exp(-\frac{1}{x})$ for $x>0$. Note that for $x\leq0$, ${\cSx}$ is  one-dimensional 
 while for $x>0$, ${\cSx}$ is a single point,  violating the definition of \flat.  That is,
${\cSx}$ does not ``have
the same dimension except for a measure zero exceptional set of $x$''. In addition, for $x<0$, ${\cSx}$ is 
$(-\infty,0]$
a ``manifold with boundary'', which also violates our definition.

The following result shows that the ACR property is a sufficient condition for flatness.

\begin{theorem}[The Flat Theorem]
\label{thm:flat}
Let $F:U\subset \mathbb{R}^N\to \mathbb{R}^M$ be a $C^\infty$ function. If $F$ is ACR, then $F$ is $d$-\flat, with dimension $d=N-$ maxrank$(F)$.
\end{theorem}

{If $F$ is ACR and $d=N-$ maxrank$(F)$, then for almost every $p\in U$, $\cSp$ is a $d$-dimensional manifold and the tangent space of $\cSp$ at $p$ is the 
kernel of $DF(x)$.}

The proof of the above theorem is given in the next section. It uses Sard's 1965 Theorem~\cite{sard_hausdorff_1965}, which is a considerable generalization of his better known 1942 result~\cite{sard_measure_1942}. 

\NI {For a function $F$ we say a point $p$ is {\bf exceptional} if rank$(DF(p)) <$ maxrank($F$). Let $E_F$ denote the set of exceptional points of $F$. The proof of Theorem~\ref{thm:flat} requires us to show not only that $E_F$ has measure zero, 
but that {the union of $\cSp$ over all $p \in E$} also has measure zero. The following example shows how $\cSp$ can be much larger than $p$ for an exceptional point $p$.}

\begin{example}\normalfont \label{ex:xy}
 A simple but nontrivial application of Theorem \ref{thm:flat} is the function $F(x_1,x_2) = x_1x_2$ from $\mathbb{R}^2$ to $\mathbb{R}^1$. The function $F$ is ACR because rank $DF(x)$ is one except at $x=(0,0)$. 
 The conclusion of the theorem is that SolSet$(p)$ is a smooth $C^\infty$ manifold of dimension one, the union of two branches of a hyperbola, except for $p$ on either of the $x_1$ or $x_2$ axes. Then, the solution set is the union of the axes, which is not a manifold. 
 
 Note that for points $p=(x_1,0)$ or $(0,x_2)$ with $x_1\neq 0$ or $x_2\neq 0$, the Implicit Function Theorem implies the manifold property {\it locally}, but the solution set SolSet(p) is globally the union of the axes -- not a manifold. Therefore these points $p$ lie in the measure-zero set where flatness fails.
\end{example}

The surprising fact about Theorem \ref{thm:flat} is that the same behaviors seen in the example follow also 
for every $C^\infty$ functions $F$, as long as the ACR property holds. Furthermore, the following result is an immediate corollary of
 Theorems~\ref{thm:key} and \ref{thm:flat}.

\begin{corollary}{}\label{cor:flat}
Assume  ${\mathcal{F}}$ is an ACR vector space of $C^\infty$ functions $F:U\to \mathbb{R}^M$. Then almost every $F\in{\mathcal{F}}$ is $d$-\flat\  with dimension $d=N-$maxrank$({\mathcal{F}})$. 
\end{corollary}
{
\section{Applications}\label{Sec:applications}
In this section we discuss examples that illustrate Corollary~\ref{cor1} and Theorem~\ref{thm:flat}.
 Systems of equations collected into vector spaces are commonplace in engineering applications. To begin, consider the following simple example of mechanical linkage problems in the spirit of those in~\cite{duka_neural_2014,siciliano_kinematic_1990,Thurston}.
\begin{example}\normalfont
\label{ex:robot}
 The robotic arm has two interior joints, at $u_1\in\mathbb{R}^3$ (the elbow) and  $u_2\in\mathbb{R}^3$ (the wrist),
 each of which is considered a variable. 
 One end is fixed at a pivot point $s_1\in\mathbb{R}^3$, the shoulder, and the other end is fixed at a pivot point $s_2\in\mathbb{R}^3$, the end of a hand.  The joint positions $u_1$ and $u_2$ satisfy the following length restrictions. 
\begin{align}
||s_1-u_1||^2 = c_1, \nonumber\\
||s_2-u_2||^2 = c_2,\\
||u_2-u_1||^2 = c_3. \nonumber
\end{align}
More generally the equations for $u_1$ and $u_2$ have the following form.
\begin{eqnarray}\label{eq:robut-arm}
f_1(u_1)&=&c_1,\nonumber\\
f_2(u_2)&=&c_2,\label{eq:f2}\\
f_3(u_1,u_2)&=&c_3.\nonumber
\end{eqnarray}
Write $p=(u_1,u_2) \in \mathbb{R}^6$  and $F=(f_1,f_2,f_3)$. 
Theorem~\ref{thm:flat} implies the following global result. Almost every $C^\infty$ function $F:\mathbb{R}^6\to\mathbb{R}^3$ that has the form (\ref{eq:robut-arm}) has the following property: For almost every $p\in\mathbb{R}^6$, SolSet$(p)$ is a manifold of dimension 3.
\label{Ex:robot}
\end{example}

\begin{example}
\rm Consider the ordinary differential equation
\begin{eqnarray*}
\dot{m}_1 &=& f_1(m_1,p_1,p_2),\\
\dot{m}_2 &=& f_2(m_2,p_1),\\
\dot{p}_1 &=& f_3(m_1,p_1),\\
\dot{p}_2 &=& f_4(m_2,p_2).
\end{eqnarray*}
modeling a two-gene regulatory network \cite{chesi2008stability}. Here $m_i, p_i$ represent the concentrations of mRNA and protein of gene $i$, for $i=1,2$. The model assumes that gene 1 is an activator/regressor of genes 1 and 2, and gene 2 is an activator/regressor of gene 1 only. To find equilibria of the network we set the left sides to zero. Let $\cF$ denote the function space of $C^\infty$ functions $F=(f_1,f_2,f_3,f_4):R^4\to R^4$ of this form, which is ACR by Proposition \ref{prop:str}. Then Corollary \ref{cor1} states that either solutions are robust for almost every $F\in \cF$, or fragile for almost every $f\in \cF$. It is easy to check that the maximum rank of the $4\times 4$ Jacobian of the system is equal to 4, so the former case holds. 
\end{example}
\begin{example} \rm
{ The following structure graph represents the ``JaK (Janus kinase
)/Stat'' signaling pathway. There are clinical evidence that confirms that Jak/Stat signaling pathway is often activated in hematologic
cancers~\cite{thomas_role_2015}. Therefore, understanding this signaling pathway could help in designing more efficient targeted therapies to suppress this pathway. 
 This graph
is motivated by Model MedB-1 in \cite{raia_dynamic_2011}. For the convenience of the reader we write the structured system associated with Model MedB-1 below.
\begin{equation}\label{eq:signal}
\begin{array}{rrcl}
{\rm Rec}  &\dot{x}_1 &=& f(x_1,x_2),\\
{\rm Rec\_i} &\dot{x}_2 &=& f_2(x_1,x_2),\\
{\rm IL13\_ Rec}&\dot{x}_3 &=& f_3(x_1,x_3,x_7),\\
{\rm p\_IL13\_Rec} &\dot{x}_4 &=& f_4(x_3,x_4,x_7),\\
{\rm p\_IL13\_Rec\_i} &\dot{x}_5 &=& f_5(x_4,x_5),\\
{\rm JAK2} &\dot{x}_6 &=& f_6(x_3,x_4,x_6,x_7,x_{11}),\\
{\rm pJAK2} &\dot{x}_7 &=& f_7(x_3,x_4,x_6,x_7,x_{11}),\\
{\rm STAT5} &\dot{x}_8 &=& f_8(x_7,x_8,x_9),\\
{\rm pSTAT5} &\dot{x}_9 &=& f_9(x_7,x_8,x_9),\\
{\rm SOCS3mRNA} &\dot{x}_{10} &=& f_{10}(x_9),\\
{\rm SOCS3} &\dot{x}_{11}&=&f_{11}(x_{10},x_{11}),\\
{\rm CD274mRNA} &\dot{x}_{12}&=&f_{12}(x_9).\\
\end{array}
\end{equation}
The above pathway starts by binding an enzyme called IL13 to receptor Rec (variable $x_1$) and ends up in the production of two mRNAs, 
CD274mRNA ($x_{12}$) and SOCS3mRNA ($x_{10})$. What triggers the transcription of these two mRNAs is a molecule called pSTAT5 ($x_9$). Once $x_9$ is produced, either it triggers the transcription of CD274mRNA or SOCS3mRNA . If a suppressor gene is used to block the transcription of CD274mRNA, then the process of translation of protein {\rm SOCS3} from {\rm SOCS3mRNA} will be a robust path and it will not be sensitive to perturbations.

{ The vector space of functions that} respects this structure has maxrank 11. Hence, {Theorem~\ref{thm:flat} implies that} almost every steady state solution to the above system lies on a one-dimensional manifold. {On the other hand,} by knocking out the node $x_{12}$ (suppressing the gene that leads to production of CD274mRNA), the above structure would be robust to perturbations. One may use such a method to assess what will happen under a certain treatment regimen. For example, blocking the production of CD274mRNA ($x_{12}$) using gene suppressors will make the alternative pathway (production of SOCS3mRNA) robust.}
 \begin{figure}
\begin{center}
 \begin{tikzpicture}
    [ultra thick]
    \node[fignew2]at(2.5cm,3.4cm)(a6){6};
    \node[fignew2]at(-4cm,0cm)(a2){2};
    \node[fignew2]at(-4cm,2cm)(a3){3};
    \node[fignew2]at(-1cm,2cm)(a4){4};
    \node[fignew2]at(1.2cm,2cm)(a5){5};
    \node[fnew]at(4cm,2cm)(a14){12};
    \node[fignew2]at(-5cm,1cm)(a1){1};
    \node[fnew]at(7cm,0.3cm)(a13){11};
    \node[fignew2]at(2.5cm,0.5cm)(a7){7};
    \node[fignew2]at(4cm,0.3cm)(a9){9};
    \node[fignew2]at(2.5cm,-1cm)(a8){8};
    \node[fnew]at(5.5cm,0.3cm)(a10){10};
    \draw[<->,>=stealth,black] (a1)-- (a2);
    \draw[->,>=stealth,black] (a1)-- (a3);
    \draw[->,>=stealth,black] (a3)-- (a4);
    \draw[->,>=stealth,black] (a4)-- (a5);
    \draw[->,>=stealth,black] (a4)-- (a6);
    \draw[<->,>=stealth,black] (a6)-- (a7);
    \draw[<->,>=stealth,black] (a4)-- (2.4,0.88);
    \draw[->,>=stealth,black] (a14)-- (a6);
    \draw[->,>=stealth,black] (a14)-- (a7);
    \draw[->,>=stealth,black] (a9)-- (a14);
    \draw[->,>=stealth,black] (a9)-- (a10);
    \draw[->,>=stealth,black] (a10)-- (a13);
    \draw[<->,>=stealth,black] (a9)-- (a8);
    \draw[->,>=stealth,black] (a7)-- (a8);
     \draw[->,>=stealth,black] (a7)-- (a9);
    \draw[->,>=stealth,black] (-3.8cm,2.2cm)-- (2.2cm,3.5cm);
    \draw[<->,>=stealth,black] (-3.8cm,1.8cm)-- (2.3cm,0.75cm);
\draw[->,>=stealth,black](a1)edge[in=-20,out=60,loop below]node[below right]{}();
\draw[->,>=stealth,black](a2)edge[in=-20,out=60,loop above]node[below right]{}();
\draw[->,>=stealth,black](a3)edge[in=-20,out=60,loop above]node[below right]{}();
\draw[->,>=stealth,black](a4)edge[in=-20,out=60,loop above]node[below right]{}();
\draw[->,>=stealth,black](a5)edge[in=-20,out=60,loop above]node[below right]{}();
\draw[->,>=stealth,black](a6)edge[in=-20,out=60,loop above]node[below right]{}();
\draw[->,>=stealth,black](a13)edge[in=-20,out=60,loop above]node[below right]{}();
\draw[->,>=stealth,black](a9)edge[in=-20,out=60,loop below]node[below right]{}();
\draw[->,>=stealth,black](a8)edge[in=-20,out=60,loop left]node[below right]{}();
\draw[->,>=stealth,black](a7)edge[in=-20,out=60,loop left]node[below right]{}();
\end{tikzpicture}
\end{center}
\caption{{\bf A signaling pathway leading to cancer.} This structure graph represents the Jack/stat signaling pathway represented by system of equations~\eqref{eq:signal}.
Every 12 by 12 matrix $S$ that respects the above structure has rank at most 11.  By Corollary~\ref{cor1} this signaling pathway is not robust.
}
\end{figure}
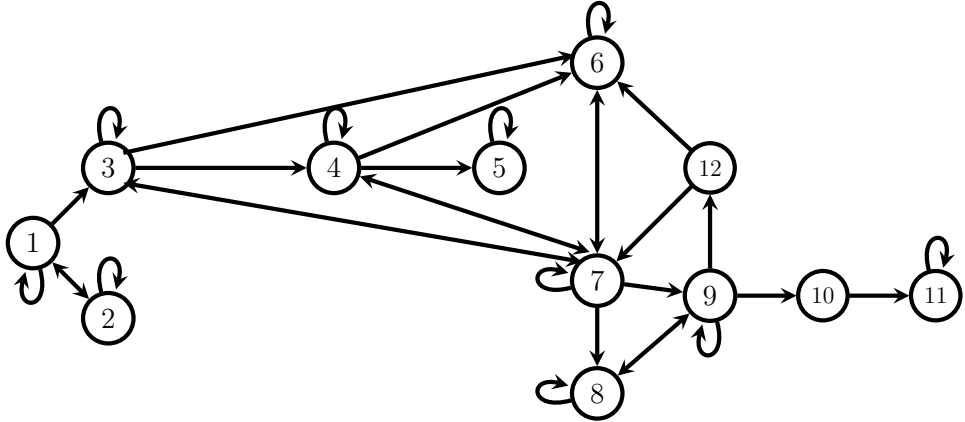
\end{example}
\begin{example}\rm
Consider the ``trophic'' ecosystem of three predator species and two prey species illustrated in Fig. \ref{fig:trophic}. This is a slightly more complicated version of Figs. 1 and 2. Let $\cF$ be the function space of structured systems of this form. Note that the Jacobian matrix has form 
 \begin{equation}
 DF = \left[
\begin{array}{ccccc}
0&0&0&f_{14}&f_{15} \\
0&0&0&f_{24}&f_{25} \\
0&0&0&f_{34}&f_{35} \\
f_{41}&f_{42}&f_{43}&0&0\\
f_{51}&f_{52}&f_{43}&0&0\\
\end{array}
\right]
 \end{equation}
which has max rank 4, and is fragile for every $F\in \cF$. According to Theorem \ref{thm:flat}, SolSet$(p)$ is a one-dimensional manifold for almost every $F\in\cF$.
\begin{figure}
\centering
 \begin{tikzpicture}
[ultra thick]
  \tikzset{every node}=[font=\bf]
  \node[regnode](c1)[xshift=-2cm]{1};
  \node[regnode](c2)[]{2};
\node[regnode](c3)[xshift=+2cm]{3};
\node[regnode](r1)[xshift=-1cm,yshift=-2cm]{4};
  \node[regnode](r2)[xshift=+1cm,yshift=-2cm]{5};
\draw[<->,>=stealth,black] (r1)-- (c1);
\draw[<->,>=stealth,black] (r1)-- (c3);
\draw[<->,>=stealth,black] (r1)-- (c2);
\draw[<->,>=stealth,black] (r2)-- (c1);
\draw[<->,>=stealth,black] (r2)-- (c3);
\draw[<->,>=stealth,black] (r2)-- (c2);

\end{tikzpicture}
    \caption{{\bf Trophic systems.} Predator species numbered 1, 2, and 3 interact with prey species 4, 5, {with no intralevel interactions}. The system represented by this graph cannot have robust equilibrium solutions. Adding one well-placed arrow, however, changes  the structured system into one with robust solutions for almost every choice of $C^\infty $ functions. \label{fig:trophic}}
\end{figure}
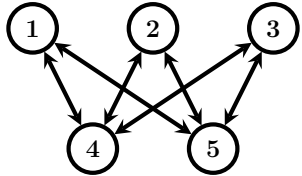

Note that adding one more connection to the network, say an arrow from predator species 1 to predator species 2, changes the max rank to 5. In this case, Corollary~\ref{cor1} implies that solutions are robust for almost every $F\in\cF$. Therefore, unlike the original system, the revised system has plausible biological solutions.
\end{example}

\begin{example}\rm
 Fig. \ref{fig:sole} shows an ecosystem proposed by Sol\'e and Montoya in \cite{sole_complexity_2001}. One can show that in this 26-species system, the max rank of the Jacobian matrix is 20. Therefore robust solutions do not exist. According to Theorem \ref{thm:flat}, when solutions do exist, they belong to 6-dimensional $C^\infty$ manifolds for almost every $F:R^{26}\to R^{26}$ in the structured function space of the system.
\end{example}

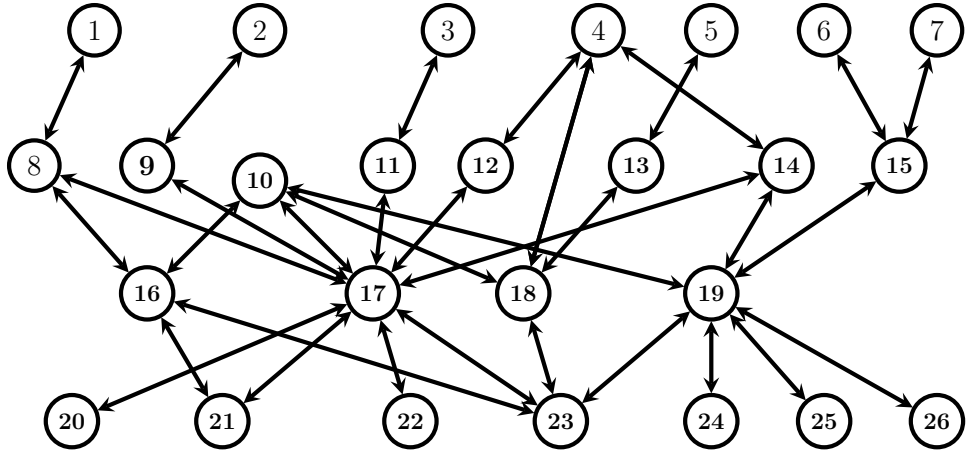
\begin{figure}
    \centering

 \begin{tikzpicture}
    [ultra thick]
     \tikzset{every node}=[font=\BF]
    \node[2uninode]at(-5.7cm,2cm)(a11){1};
    \node[2uninode]at(-3.5cm,2cm)(a12){2};
    \node[2uninode]at(-1cm,2cm)(a13){3};
    \node[2uninode]at(1cm,2cm)(a14){4};
    \node[2uninode]at(2.5cm,2cm)(a15){5};
    \node[2uninode]at(4cm,2cm)(a16){6};
    \node[2uninode]at(5.5cm,2cm)(a17){7};
    \node[2uninode]at(-6.5cm,.2cm)(a21){8};
    \node[2uninode]at(-5cm,.2cm)(a22){${9}$};
    \node[uninode]at(-3.5cm,0cm)(a23){${10}$};
    \node[uninode]at(-1.8cm,.2cm)(a24){${11}$};
    \node[uninode]at(-0.5cm,.2cm)(a25){${12}$};
    \node[uninode]at(1.5cm,.2cm)(a26){${13}$};
    \node[uninode]at(3.5cm,.2cm)(a27){${14}$};
    \node[uninode]at(5cm,.2cm)(a28){${15}$};
     \node[uninode]at(-5cm,-1.5cm)(a31){${16}$};
    \node[uninode]at(-2cm,-1.5cm)(a32){${17}$};
    \node[uninode]at(0cm,-1.5cm)(a33){${18}$};
    \node[uninode]at(2.5cm,-1.5cm)(a34){${19}$};
    \node[uninode]at(-6cm,-3.2cm)(a41){${20}$};
    \node[uninode]at(-4cm,-3.2cm)(a42){${21}$};
    \node[uninode]at(-1.5cm,-3.2cm)(a43){${22}$};
    \node[uninode]at(0.5cm,-3.2cm)(a44){${23}$};
    \node[uninode]at(2.5cm,-3.2cm)(a45){${24}$};
    \node[uninode]at(4cm,-3.2cm)(a46){${25}$};
    \node[uninode]at(5.5cm,-3.2cm)(a47){${26}$};
  \draw[<->,>=stealth,black] (a11)-- (a21);
    \draw[<->,>=stealth,black] (a12)-- (a22);
    \draw[<->,>=stealth,black] (a13)-- (a24);
    \draw[<->,>=stealth,black] (a14)-- (a25);
    \draw[<->,>=stealth,black] (a14)-- (a27);
    \draw[<->,>=stealth,black] (a14)-- (a33);
    \draw[<->,>=stealth,black] (a14)-- (a33);
    \draw[<->,>=stealth,black] (a15)-- (a26);
    \draw[<->,>=stealth,black] (a16)-- (a28);
    \draw[<->,>=stealth,black] (a17)-- (a28);
    \draw[<->,>=stealth,black] (a21)-- (a31);
    \draw[<->,>=stealth,black] (a21)-- (a32);
     \draw[<->,>=stealth,black] (a22)-- (a32);
    \draw[<->,>=stealth,black] (a23)-- (a31);
     \draw[<->,>=stealth,black] (a23)-- (a32);
     \draw[<->,>=stealth,black] (a23)-- (a33);
     \draw[<->,>=stealth,black] (a23)-- (a34); 
    \draw[<->,>=stealth,black] (a24)-- (a32);
    \draw[<->,>=stealth,black] (a25)-- (a32);
    \draw[<->,>=stealth,black] (a26)-- (a33);
    \draw[<->,>=stealth,black] (a27)-- (a32);
 \draw[<->,>=stealth,black] (a27)-- (a34);
  \draw[<->,>=stealth,black] (a28)-- (a34);
  \draw[<->,>=stealth,black] (a34)-- (a47);
  \draw[<->,>=stealth,black] (a34)-- (a46);
  \draw[<->,>=stealth,black] (a34)-- (a45);
  \draw[<->,>=stealth,black] (a34)-- (a44);
  \draw[<->,>=stealth,black] (a33)-- (a44);
  \draw[<->,>=stealth,black] (a31)-- (a44);
  \draw[<->,>=stealth,black] (a31)-- (a42);
  \draw[<->,>=stealth,black] (a32)-- (a42);
  \draw[<->,>=stealth,black] (a32)-- (a41);
  \draw[<->,>=stealth,black] (a32)-- (a43);
  \draw[<->,>=stealth,black] (a32)-- (a44);
\end{tikzpicture}

\caption{The graph represents a proposed ecosystem of Sol\'e and Montoya in \cite{sole_complexity_2001}. The max rank of the Jacobian of the equilibrium equations for the 26 species is 20. No solutions are robust. By Theorem \ref{thm:flat}, for almost every $F$ of this form, for almost every $p$, solution set solset$(p)$ is 6-dimensional $C^\infty$ manifolds in the domain $R^{26}$}
\label{fig:sole}
\end{figure}
}

 \section{Proofs of Theorems} \label{s:proofs}
 In this section, we prove { Theorems~\ref{thm:key} and \ref{thm:flat} }. We will find versions of Fubini's Theorem helpful in several ways.

To make clear what ``almost every'' and ``measure 0'' mean, we sometimes write ${\leb}$-almost every when we mean Lebesgue measure almost every, and write ${\leb}_d$ to denote $d$-dimensional Lebesgue measure. We denote  Lebesgue measure on a finite dimensional Euclidean space $Y$ by ${\leb}(Y)$.

\NI
{\bf Fubini Theorem.} Let $X$ and $Y$ be finite dimensional Euclidean spaces. Let $E$ be a {measurable} subset of $X\times Y$.  Then $E$ has ${\leb}(X\times Y)$ measure zero if and only if
 for ${\leb}(Y)$ almost every $y\in Y$, $(X\times \{y\}) \cap E$ has ${\leb}(X)$ measure zero.

\NI 
In the following result, there are infinitely many linear coordinate choices one could use to define Lebesgue measure on $Y$, but the resulting Lebesgue measures all agree on which sets have measure 0.

\NI
{\bf Fubini Corollary.} Let $X$ and $Y$ be finite dimensional vector spaces with {$X \subseteq Y$}. 
 Let $E$ be a measurable subset of $Y$.\\
Then $E$ has measure zero in $Y$ if and only if
\begin{equation}\label{fubini X}
(X+y) \cap E\text{ has }  \leb (X)\text{ measure zero for } \text{ almost every } y\in Y,
\end{equation}
where  $\leb (X)$-measure zero means
with respect to the translated Lebesgue measure on $X+y$.

\NI
 We will use the above Fubini Corollary by showing that if there is a subspace $X$ for which 
the property (\ref{fubini X}) holds, then it holds for every subspace $X$.
Below in the proof of Theorem~\ref{thm:key}, we apply this indirectly for $Y$ of dimension $N+1$ and we use two choices of $X$, one with dimension 1 and the other with dimension $N$.

The following maximal rank result for a vector space of linear functions is elementary and is included as 
an illustrative example of an almost every property and as a method of proof we use later.
 In particular,  the vector space $V$ in Proposition \ref{prop:matrices} is ACR.
 
 \begin{proposition}
 \label{prop:matrices}
 Let $V$ be a subspace of $M\times N$ matrices $\mathbb{ R}^{MN}$. Let $\rho = \max_{A\in V}{ \rm rank}(A)$.  Then almost every matrix $A$ in $V$ has rank $\rho$ (``almost every'' with respect to Lebesgue measure on $V$).  
 \end{proposition} 
 
\begin{proof}
 Let $B$ be a matrix in $V$ with maximal rank$(B)=\rho$. We will show that for each $A\in V$, rank$(A+cB) =\rho $ for almost every scalar $c$. 
There is at least one nonsingular $\rho\times\rho $ submatrix $B_\rho$ of $B$. The corresponding $\rho\times\rho $ submatrix of $A+cB$ has determinant that is a degree $\rho$ polynomial in $c$. We need to show the polynomial is not identically zero.
It equals $c^\rho\cdot\det_\rho 
 (\frac{1}{c}A_\rho+ B_\rho)$ for $c\ne0$.
 It is not identically zero since for large $c$,
 $\det
 (\frac{1}{c}A_\rho+ B_\rho)$  approaches $\det_\rho(B_\rho)$, which by assumption is nonzero.
By the Fundamental Theorem of Algebra, $A+cB$ has rank $\rho$ for all but a finite set of $c$.  
Now the Fubini Corollary applies where $y=A$, $Y=V$, and $X$ is the  one-dimensional subspace including $B$, and $E$ is the exceptional subset of $V$ of matrices with rank $<\rho$. \end{proof}

\begin{proof}[Proof of Theorem~\ref{thm:key}]
Let $\rho=\text{maxrank}({\mathcal{F}})$. Let $\mF$ be a rank-maximizer function in ${\mathcal{F}}$. Then, rank$(D\mF(x)) =\rho$ for almost every $x\in U$. 

Let Fubini Corollary subspace $X$ be the one-dimensional subspace of ${\mathcal{F}}$ consisting of the functions $c\mF$ for $c\in\mathbb{R}$.
Define $U_\rho := \{x: \rank(D\mF(x))=\rho\}$.
 To prove the lemma, it is sufficient to show that for each $F\in{\mathcal{F}}$,
$F+c\mF$ is a rank-maximizer for almost every $c\in\mathbb{R},\ie,$
for each $F$, for almost every $x \in U_\rho$,
\begin{equation}\label{claim c}
 \rank \qty( D(F+c\mF)(x))=\rho \text{ for  almost every }  c\in{\mathbb{R}}.
\end{equation}
Fix an $F$.  Let $E=\{(x,c):\rank( D(F+c\mF)(x)) <\rho\}$ and\\
$E^x=\{(y,c)\in E:y=x\}$.

Since $c \in \mathbb{R}$, we will first show $E^x$ has one-dimensional measure zero for each $x\in U_\rho$.
For each $x\in U_\rho$, apply the argument in the proof of Proposition~\ref{prop:matrices} (where a $\rho \times\rho $ submatrix is chosen) and conclude that
 $E^x$ contains at most finitely many points and so has measure 0.
 Apply the Fubini corollary to conclude
that since {Lebesgue} almost every ``slice'' $E^x$ of $E$ has measure 0, so does $E$.

Next we switch the roles of $x$ and $c$ and prove that 
for almost every $c$, 
\begin{equation}\label{claim x}
\rank\big(D(F+c\mF)(x)\big)=\rho\text{ for  almost every } x\in U.
\end{equation}
That is, $F+c\mF$ is a rank-maximizer for almost every $c\in\mathbb{R}$.
 
Define
$E_c:=\{(x,s)\in E:s=c \}$. 
 Above by using Fubini's Theorem we showed that the set $E$ has measure zero, which implies by Fubini's Theorem that for almost every $c\in\mathbb{R}$, $E_c$  has measure zero, proving \eqref{claim x}.
\end{proof}

The Hausdorff dimension of a set can be defined by first defining measure zero for a Hausdorff measure of dimension $s$, which is an extension of Lebesgue measure zero. 
\begin{definition}\rm 
We say that $B$ has {\BF \bf $s$-Hausdorff measure zero} if for each $\varepsilon > 0$ there is a countable cover of $B$ by sets for which the sum of the $s$-th powers of the ``diameters'' of the covering sets is less than $\varepsilon$. The diameter of a set $C$ is $\sup\|x-y\|$ over $x,y\in C$.
Let $d$ be the infimum of $s$ for which $B$ is $s$-Hausdorff measure zero. Then $d$ is the {\bf Hausdorff dimension} of $B$.  
{We will also encounter $\rho$-dimensional Hausdorff measure on a manifold of dimension $\rho$, in which case that Hausdorff measure equals Lebesgue measure of the manifold.}
\end{definition}

{\bf Notation for Lemma~\ref{lemma measure} and its proof.}  
Let $F:U\subset\mathbb{R}^N\to\mathbb{R}^M$ be a $C^{\infty}$ function, and
define $\rho:=\displaystyle{\max_{x\in U}\ \rank(DF(x))}$.
Define  $U_\rho:= \{x \in U: \rank(DF(x))=\rho \}$, an open subset of  $\mathbb{R}^N$.
As mentioned above, ${\leb}_k$ denotes $k$-dimensional Lebesgue measure.

Define $A_{\rho-1}:=\{F(x): x\in U \text{ and } \rank(DF(x))\leq\rho-1 \}$.  Sard \cite{sard_hausdorff_1965} proved that $A_{\rho-1}$ has Hausdorff dimension $\le\rho -1$. Of course there may also be $x'$ for which $F(x')\in A_{\rho-1}$ and $\rank(DF(x'))=\rho$.

\begin{figure}[H]
 \centering
\begin{tikzpicture}
[scale=0.5]
\draw [white,fill = red!40, fill opacity=.3](-14,2) -- (-14,5) -- (-11,5) -- (-11,2)-- (-14,2)  -- cycle;
\node[]at(-12.2,4.5){$B$};
\draw[->] (-16,0) -- node[pos=.95,below]{$\mathbb{R}^{N-\rho}$}(-8,0);
\draw[->] (-16,0) -- node[pos=.95,left]{$\mathbb{R}^\rho$}(-16,7);
\draw[->] (-4,0) -- node[pos=.95,below]{$\mathbb{R}^{M-\rho}$}(4,0);
\draw[->] (-4,0) -- node[pos=.95,left]{$\mathbb{R}^\rho$}(-4,7);
\draw[-,line width=2,green](-15.7,3.1) --(-8.5,3.1);
\node[]at(-9.5,3.5){$x+Y$};
\draw[-,line width=2,red](-16,2) -- node[pos=.65,left,black]{$B^*$}(-16,5);
\draw[-,line width=2,red](-4,2) -- node[pos=.65,right,black]{$F(B^*)=F(B)$}(-4,5);
\end{tikzpicture}
\caption{{\bf Sketch of Lemma~\ref{lemma measure}.} Projection of the set $B$  where $F$ is defined in (\ref{projection}).} \label{fig:projection}
\end{figure}
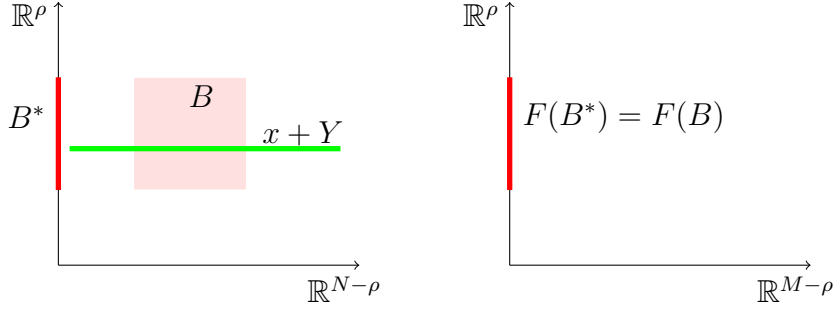

For any set $A\subset\mathbb{R}^M$, let $F^{-1}(A)=\{x:F(x)\in A\}$.
 For the proof of Theorem~\ref{thm:flat} we need to show that $U_\rho\cap F^{-1}(A_{\rho-1})$ has measure 0 in $\mathbb{R}^N$.
\begin{lemma}\label{lemma measure}
Let $A\subset\mathbb{R}^M$ have 
$\rho$-Hausdorff measure 0.
Then $U_\rho\cap F^{-1}(A)$ has Lebesgue measure 0.
\end{lemma}
\begin{proof}
Assume the notation in the Lemma.
We will prove that if 
$A$ has $\rho$-Hausdorff measure 0
 and $B= U_\rho \cap F^{-1}(A)$, then $B$ has measure 0 in $\mathbb{R}^N$. 

{ Suppose the contrary, that $B$ has positive ${\leb}_N$-measure. 
 Choose a point $q \in U_\rho$ so that every neighborhood $U_q$ of the point 
$p$ intersects $B$ in a set of positive ${\leb}_N$-measure. 
The Constant Rank Theorem (Theorem 11.1 in \cite{Loring}) says that under the above hypotheses, there is a neighborhood $U_q$ of $q$ in $U_\rho$ on which $F$ is the projection:
\begin{align}\label{projection}
F(x_1,\ldots,x_N)= (x_1,\ldots,x_\rho ,0,\dots,0)\in\mathbb{R}^M
\end{align}
}
for some smooth choice of coordinates in the domain $U_q$ and range $\mathbb{R}^M$. 
We apply the Fubini Theorem to the projection (\ref{projection}).
Let $X$ be the subspace of $\mathbb{R}^N$ of points $(x_1,\ldots,x_\rho,0,\dots,0)$ and $Y$ the complementary space of points $y=(0,\dots,0,y_{\rho+1},\ldots,y_N)$.
Let $B^*=\{x \in X:(x+Y)\cap (B \cap U_q) \text{ has positive ${\leb}_{N-\rho}$-measure}\}$. 
Fubini says $B$ having positive ${\leb}_{N}$-measure, implies
$B^*$
has positive ${\leb}_{\rho}$-measure in $Y$, as sketched in Fig.\ref{fig:projection}. 
 
 The $\rho$-dimensional Hausdorff measure of a set $F(B\cup U_q)$ is its $\rho$-dimensional Lebesgue measure on the subspace $Z_\rho=\{(z_1,\ldots,z_\rho ,0,\dots,0):\text{for all }z_1,\ldots,z_\rho\in\mathbb{R}\}\subset\mathbb{R}^M$, 
 and $F(B^*)$ has the same measure ($\rho$ -dimensional Lebesgue measure) as $B^*$. 
 This is the $\rho$-dimensional Hausdorff measure of $F(B\cup U_q)$ in $\mathbb{R}^M$.
 We conclude that
 $F(B^*)$ has positive measure, and therefore $A$ has positive measure since
 $F(B)\subset A$, contradicting our assumption that $A$ has measure zero.
 \end{proof}
 
\begin{lemma}\label{prop:manifold} 
Let $F:U\subset\mathbb{R}^N\to\mathbb{R}^M$ be a $C^\infty$ function.
 Let $\rho= \text{maxrank}(F):={\max_{x\in U}\rank(DF(x))}$.
 Assume some $p\in U$ has the following property: $DF(x)$ has rank $\rho$
for all $x\in\cSp$. Then $\cSp$ is a manifold of dimension $N-\rho$.
\end{lemma}

\begin{proof}
Note that if $DF(p)$ has rank $\rho$, then there are $\rho$ vectors in $\mathbb{R}^N$ whose images under $DF(p)$ are linearly independent. That is an ``open'' property in the sense that $p$ has an open neighborhood in which those vectors are independent. 
Since $\rho$ is the maximum of rank$(DF(x))$, there is an open  neighborhood of $\cSp$ on which the 
the rank of $DF(x)$ is $\rho$. 
This is the precise setting of Thm.~11.2 of the ``Constant-rank level set theorem'' in~\cite{Loring}. There $F$ is assumed to be a $C^\infty$ map between manifolds, which in our case are $U$ and $\mathbb{R}^M$.
\end{proof}
 Next we give a proof of Theorem~\ref{thm:flat}.
When $\rho=M$, the proof can be simplified and follows from Sard's 1942 theorem \cite{sard_measure_1942}. The general case, including $\rho<M$, requires Sard's 1965 theorem~\cite{sard_hausdorff_1965} involving Hausdorff dimension.

 \begin{proof}[{\it Proof of Theorem~\ref{thm:flat}}]
Assume $F$ is $C^\infty$ and is ACR.
Let $\rho = $maxrank$(F)$.
Following Sard~\cite{sard_hausdorff_1965}, for each integer $n\ge0$, we define the set 
$$A_n=\{c\in\mathbb{R}^M: \text{ there exists an } x \text{ such that } F(x)=c \text{ and } \rank\big(DF(x)\big)\le n\}.$$ 
When $\rho$ is the largest value of $\rank(DF)$, as it is in our case, {then the $\rho$-Hausdorff measure of  $A_\rho$ is positive, 
and by Sard's 1965 paper, $A_\rho$ has dimension $\rho$ while 
$A_{\rho-1}$ has dimension at most $\rho-1$, and its $\rho$-Hausdorff measure is 0.} 
We need $F^{-1}(A_{\rho-1})\cap U_\rho$ to have measure 0 in $\mathbb{R}^N$. 
This is assured by Lemma~\ref{lemma measure}. 
Hence for almost every point $p\in U_\rho$, $\cSp \subset U_\rho$, in which case all points $q \in \cSp$ satisfy 
$\rank(DF(q))=\rho$. 
By Lemma~\ref{prop:manifold},
$\cSp$ is a manifold of dimension  $N-\rho$. \end{proof}

\section{Discussion}\label{s:discussion}

{Structured systems of equations are common in scientific applications where detailed information about the relationships between variables is scarce. Section \ref{Sec:applications} shows several examples of structured systems in ecology, genetic networks, and engineering systems. 
}







Our interest {in this study was to explore the implications of such  structure on } generic properties of solutions, regardless of the specific functions involved, {in the same manner that the Competitive Exclusion Principle imposes limits on what solutions can be robust}. By generic properties, we mean properties 
held by almost every function in the vector space of functions with that structure.

{To better understand global properties, we coined the term ``almost constant rank'' (ACR).
It turned out that this {property, possessed by all linear systems,} was the pivotal property that {encodes the common properties that} vector spaces of structured systems have.
 Our main results are stated in the more general context of ACR vector spaces of functions, of which structured systems are examples.}

The theory presented here is related to and can be compared with Sard's 1942 and 1965 Theorems and the Implicit Function Theorem.

\NI {\bf Comparison with Sard's Theorems.}
 Let $U\subset \mathbb{R}^N$ be open and 
let $F:U\to\mathbb{R}^M$.
Theorem~\ref{thm:flat} can be viewed as a dual of (S2), below, which is a consequence of 
Sard’s 1965 Theorem~\cite{sard_hausdorff_1965}. Our statement is about the domain $U\subset \mathbb{R}^N$ of a function $F$, while Sard's is about the range space $\mathbb{R}^M$.

Let $Z_F:=\{x\in U:\rank(DF(x))< M\}$. 
Theorem~\ref{thm:key} concludes under certain circumstances that for almost every $F$,
$Z_F$ has Lebesgue measure zero. In contrast,
Sard's 1942 theorem, which we denote by $S_{1942}$,
says its image $F(Z_F)$ has Lebesgue measure zero~\cite{sard_measure_1942}.

\NI
Write $F^{-1}(c)=\{x:F(x)=c\}$. Sard's  $S_{1942}$  together with our Prop~\ref{prop:manifold} implies the following.

\NI
{\bf \BF $(S1)=(S_{1942} +$ Prop.~\ref{prop:manifold})} 
For (Lebesgue) almost every  $c\in\mathbb{R}^M$, $F^{-1}(c)$ either is a manifold of dimension $N-M$ or is the empty set.

\NI
Let $\rho  = $ max$_x$ rank$(DF(x))$. When $\rho<M$, (S1) tells us nothing about the solutions. Since in this case, for almost every $c$, $F^{-1}(c)$ is empty. 
In 1965 Sard generalized his 1942 result in a manner that is important for us.
Sard's 1965 theorem, which we denote by $S_{1965}$, says the following:
{ Let $X_{\rho-1} = \{x\in U\subset\mathbb{R}^N: DF(x) \text{ has rank}\leq\rho - 1\}$. Then the set $F(U)$ has Hausdorff dimension $\rho$ and $F(X_{\rho-1})$ has Hausdorff dimension $\le\rho -1.$} 
That implies that
almost every $c\in F(U)$ (with respect to 
$\rho$-dimensional Hausdorff measure), $c\notin F(X_{\rho-1})$. The following is a consequence.

As argued in \cite{sard_hausdorff_1965},
since max$_x$ rank$(DF(x))=\rho$, $\rho$ is the Hausdorff dimension of $F(U)$,  and the $\rho$-dimensional Hausdorff measure of $F(U)$ is positive and possibly infinite. 

\NI
{\bf\BF 
$(S2)=(S_{1965} +$ Prop.~\ref{prop:manifold})
} 
For almost every $c \in F(U)$ (``almost every'' with respect to $\rho$-dimensional Hausdorff measure), whenever $F(x)=c$, $DF(x)$ has rank $\rho$. Therefore, $F^{-1}(c)$ is a manifold of dimension $N-\rho$.

{

\NI
{\bf The role of the Implicit Function Theorem (IFT).}
Theorem \ref{thm:flat} can be viewed as a globalization of the IFT. If it is known that the Jacobian $DF(p)$ maps onto the target tangent space for almost every $p$ in the domain, then the IFT shows that locally, the solution set has smooth manifold structure in a neighborhood of such points $p$. 
However, as Example \ref{ex:xy} shows, this does not mean that the solution set is a manifold, globally speaking. The problem is that as the solution set is followed beyond the local neighborhood, a point may be encountered in the solution set where the Jacobian rank drops.  
The fact that this is almost always avoided is precisely the extra information that Theorem \ref{thm:flat} provides. 
When combined with Theorem \ref{thm:key}, we find that this behavior is actually prevalent, in the formal sense, in vector spaces with the ACR property.

\NI
In this treatment, we have restricted discussion to vector spaces of functions $F:U\subset \mathbb{R}^N \to \mathbb{R}^M$, for simplicity. However, the proofs extend almost without change to vector spaces of functions between $C^\infty$-manifolds of dimensions $N$ and $M$.
See~\cite{Arnold1971} for a parametric approach to the problem of rank.}

The main theoretical results of this article, Corollary~\ref{cor1} and Theorem \ref{thm:flat}, show that generic properties of solutions of structured systems depend crucially on a single number, the generic rank of the structure matrix $S$ of the system. The rank is in turn connected to the topological properties of the associated directed graph of the system, which will be addressed more fully in future work.

\section*{Acknowledgments}
We thank Dima Dolgopyat, Shuddho Das, and Roberto De Leo for their helpful comments.

\appendix 
\section{A brief overview of prevalence}\label{prevalence}

\NI
The concept of prevalence  is useful when a vector space $\mathcal{F}$ is infinite dimensional.
Prevalence is a concept that is used to extend the idea of ``Lebesgue almost every'' to infinite-dimensional vector spaces~\cite{sauer_embedology_1991,hunt_prevalence_1992,ott_prevalence_2005}. 
The term ``prevalence'' was introduced by Sauer, Yorke, and Casdagli~\cite{sauer_embedology_1991} and generalized in \cite{hunt_prevalence_1992}. For a 1972 similar definition by Christensen see also~\cite{Christensen,hunt1993}. 
{Prevalence can be compared with full measure in finite dimensions from the Fubini Corollary in Section~\ref{s:proofs}. There we stated the Corollary for determining if a set was measure 0. Here we state it for full measure sets. The wording is chosen so that the prevalence definition is a small change in the wording of the Corollary, following the old practice of turning a theorem (or corollary) into a definition.

\NI {{\bf Fubini Corollary 2.} Let $Y$ and $X\subset Y$ be vector spaces
where $X$ and $Y$ are finite dimensional. Let $G$ be a measurable subset of $ Y$.\\ 
Then the corollary says ``almost every $y\in Y$ is in $G$'' \\
if
for  almost every $p\in Y$ and almost every $x\in X$, $x+p$ is in $G$.}

Infinite dimensional spaces $Y$ have no Lebesgue measure to give meaning to ``{\it almost every} $p\in Y$'' so we substitute ``{\it every} $p\in Y$''. We can extend the definition of measurable by saying $G\subset Y$ is  measurable in the sense of prevalence if for each finite-dimensional plane $Z\subset Y$, $Z\cap G$ is measurable.

\NI {\bf Prevalence Definition.} Let $Y$ and $X\subset Y$ be vector spaces where $X$ is finite dimensional. Let $G$ be a measurable subset of $Y$.\\
Then we define
``almost every $y\in Y$ is in $G$'' (in the sense of prevalence)\\
if for every $p\in Y$ and almost every $x\in X$, $x+p$ is in $G$.
}

\NI The following are some examples of how prevalence can be helpful in getting a general insight about common properties in an infinite dimensional vector space.

\NI (P1) For $1 < p \le\infty$ 
almost every sequence $(a_i)\in l^p({\mathbb{R}})$ has the property that $\sum_{i=1}^\infty a_i$ diverges. Here, a convenient probe space is the one-dimensional subspace spanned by the infinite sequence of  $1$'s.

\NI (P2) In 1994 Hunt~\cite{Hunt_1994} showed that almost every (continuous) function in $C([0,1])$ is nowhere differentiable. The proof requires a two-dimensional probe space.

\NI (P3) The following is a generalization of the Whitney Embedding Theorem.
Whitney proved that if $Q$ is a manifold of dimension $d$ and $M>2d$, then 
for topologically generic $F:\mathbb{R}^N\to\mathbb{R}^M$, $F$ is one-to-one on $Q$. That has a prevalence generalization that is useful in investigating chaotic attractors. Whitney's manifold is replaced by an arbitrary compact set $Q\subset\mathbb{R}^N$ of box dimension $d$.  Then  for ``almost every'' $F:\mathbb{R}^N\to\mathbb{R}^M$,  $F$ is one-to-one on $Q$.  In particular, the dimension $d$ does not need to be an integer~\cite{sauer_embedology_1991}.
There are many other examples of application of prevalence in~\cite{sauer_embedology_1991,hunt_prevalence_1992,hunt1993,ott_prevalence_2005}.

\NI Prevalence has two important properties that make it a useful extension of probabilistic almost every: 
(1) If the vector space $V$ is finite-dimensional, prevalence is the same as ``almost every in the sense of Lebesgue measure'', and (2) prevalence implies dense, so that it meets and exceeds a common topological version of typical. 
Even in finite-dimensional Euclidean spaces, dense (and residual) sets can have arbitrarily small measure, so the stronger property of prevalence is a handy addition to characterize typical behavior in infinite-dimensional spaces. 

\bigskip \noindent
{\bf Acknowledgment.} We thank Dima Dolgopyat, Shuddho Das, and Roberto De Leo for their helpful comments.

\end{document}